\documentclass{amsart}
\usepackage{amsmath}
\usepackage{paralist}
\usepackage{amsfonts}
\usepackage{amssymb}
\usepackage{amsthm}
\usepackage{amscd}
\usepackage{amsrefs}
\usepackage{float}
\usepackage{tikz}
\usepackage{graphicx}
\usepackage[colorlinks=true]{hyperref}
\hypersetup{urlcolor=blue, citecolor=red}
\usepackage{hyperref}

\newtheorem{theorem}{Theorem}[section]

\newtheorem{proposition}{Proposition}

\theoremstyle{definition}

\title[Sharp Strichartz estimates in spherical coordinates] 
      {Sharp Strichartz estimates in spherical coordinates}

\author[R. Schippa]{Robert Schippa}
\address{Fakult\"at f\"ur Mathematik, Universit\"at Bielefeld, Postfach 10 01 31, 33501 Bielefeld, Germany}
\email{robert.schippa@uni-bielefeld.de}

\begin{document}
\maketitle

\bigskip

\begin{abstract}
We prove almost sharp Strichartz estimates found after adding regularity in the spherical coordinates for Schr\"odinger-like equations relying on estimates involving spherical averages. Sharpness is discussed making use of a modified Knapp-type example.
\end{abstract}

\section{Introduction}
Since we consider homogeneous estimates, we shall confine ourselves to homogeneous equations:
 \begin{equation}
 \label{eq:homogeneousEquation}
 \left\{\begin{array}{cl}
 i \partial_t u(t,x) + \varphi(D) u(t,x) &= 0, \, (t,x) \in \mathbb{R} \times \mathbb{R}^n, D=(-\Delta)^{1/2}, \\
 u(0,\cdot) &= u_0.  \end{array} \right.
 \end{equation}
In the following we will deal with Schr\"odinger-like equations for most of the time, that is the dispersion relation $\varphi \in C^{\infty}((0,\infty), \mathbb{R})$ is given by
\begin{equation}
\label{eq:schroedingerEquation}
\varphi(\rho) = \rho^a, \; a>1.
\end{equation}
Strichartz estimates capture the dispersive properties of solutions to linear dispersive equations and  classical homogeneous estimates for \eqref{eq:homogeneousEquation} in the case of \eqref{eq:schroedingerEquation} state as follows
 \begin{equation}
 \label{eq:homOrdStrEst}
 \Vert u \Vert_{L_t^q L_x^p(\mathbb{R}^n)} \lesssim_{a,n,p,q} \Vert u_0 \Vert_{\dot{H}^{s}},
 \end{equation}
 where the derivatives are determined by scaling
 \begin{equation}
 \label{eq:conditionsons} 
 s=n \left(\frac{1}{2} - \frac{1}{p}\right) - \frac{a}{q}.
 \end{equation}
Keel and Tao proved the sharp range of homogeneous estimates in \cite{TaoKeel1998}, which is given by
\begin{equation}
\label{eq:ordRange}
 \frac{1}{q} \leq \frac{n}{2} \left( \frac{1}{2} - \frac{1}{p} \right), \; q,p \geq 2, \; p \neq \infty.
\end{equation}
Making use of the bilinear interpolation argument worked out in \cite{TaoKeel1998}, Cho, Ozawa and Xia showed homogeneous estimates for more general dispersion relations in \cite[Theorem~2,~p.~1123]{ChoOzawaXia2011}.\\
Sharpness of \eqref{eq:ordRange} is seen considering a Knapp-type example, that is a maximally anisotropic propagating wave. This example was already considered by Strichartz in his seminal paper \cite[Lemma~3.,~p.~707]{Strichartz1977}, in which he linked special cases of estimates of the kind \eqref{eq:homOrdStrEst} to Fourier restriction estimates through duality. We review a modified example in Proposition \ref{prop:sharpRegularity}.\\
This raised the question, whether one can extend the range of integrability coefficients found in \cite{TaoKeel1998} and \cite{ChoOzawaXia2011} if one punishes anisotropic propagation by considering angular regularity of the initial data. More precisely, for equations \eqref{eq:homogeneousEquation} in the case of \eqref{eq:schroedingerEquation} we want to consider estimates of the kind:
\begin{equation}
\label{eq:homRegStrEst}
\Vert u \Vert_{L_t^q L_x^p} \lesssim_{a,n,p,q,\alpha} \Vert u_0 \Vert_{\dot{H}^{s,\alpha}_{\omega}},
\end{equation}
where the Sobolev spaces $\dot{H}^{s,\alpha}_\omega = D^{-s} \Lambda_\omega^{-\alpha} L^2(\mathbb{R}^n) $ with angular regularity $\alpha$ are defined making use of the inhomogeneous Laplace-Beltrami operator $\Lambda_\omega = (1-\Delta_\omega)^{1/2}$, when
$$\Delta_\omega = \sum_{1 \leq i < j \leq n} \Omega_{ij}^2, \; \Omega_{ij} = x_i \partial_j - x_j \partial_i. $$
This is in fact the case and Cho and Lee found the following theorem to hold:
 \begin{theorem}[{\cite[Theorem~1.2.,~p.~994]{Cho2013}}]
 \label{thm:radStrEstA}
 Let $n \geq 2$ and suppose that $(q,p)$ satisfies 
\begin{equation}
\label{eq:extRange}
\frac{n}{2} \left( \frac{1}{2} - \frac{1}{p} \right) < \frac{1}{q} \leq \frac{2n-1}{2} \left( \frac{1}{2} - \frac{1}{p} \right), \; q,p \geq 2,
\end{equation} 
$(n,q,p) \neq (2,2,\infty) \mbox{ and } (q,p) \neq (2,(4n-2)/(2n-3))$.\\
Then we find the estimate \eqref{eq:homRegStrEst} to hold for the solution $u$ to \eqref{eq:homogeneousEquation} with $s$ from \eqref{eq:conditionsons} provided that $\alpha > ((2n-1)/(2n-2))(2/q+n/p-n/2)$.
 \end{theorem}
Sharpness of the $(q,p)$-range (also for more general dispersion relations) up to endpoints was also seen in \cite[Section~4.1.,~pp.~1005f]{Cho2013}.\\ 
  Probing the angular regularity with a Knapp-type example we find also for the more general dispersion relations that an angular regularity $\alpha = \frac{2}{q} + \frac{n}{p} - \frac{n}{2}$ is necessary:
  \begin{proposition}
  \label{prop:sharpRegularity}
  Suppose that $n \geq 2$ and $(q,p)$ satisfies \eqref{eq:extRange}. Then we find that $\alpha \geq \frac{2}{q} + \frac{n}{p} - \frac{n}{2}$ is necessary for estimate \eqref{eq:homRegStrEst} to hold.
  \end{proposition}
 Making use of a result due to Guo from \cite{Guo2014} we prove the following result establishing estimates of the kind \eqref{eq:homRegStrEst} for Schr\"odinger-like equations with sharp angular regularity up to endpoints:
 \begin{theorem}
 \label{thm:sharpRegularity}
 Let $q,p \geq 2$, suppose that we are in the case of \eqref{eq:schroedingerEquation} and
 \begin{equation*}
 \label{eq:conditionsOnQRadStrEstC}
 \begin{split}
 &\frac{n}{2} \left( \frac{1}{2} - \frac{1}{p} \right) < \frac{1}{q} < \frac{2n-1}{2} \left( \frac{1}{2} - \frac{1}{p} \right) \mbox{ for } n =2, \\
 &\frac{n}{2} \left( \frac{1}{2} - \frac{1}{p} \right) < \frac{1}{q} \leq \frac{2n-1}{2} \left( \frac{1}{2} - \frac{1}{p} \right), \; (q,p) \neq \left( 2,\frac{4n-2}{2n-3} \right) \mbox{ for } n > 2.
 \end{split}
 \end{equation*}
 Then, for $n=2$ or $n>2$ and $q = 2$, we find the estimate \eqref{eq:homRegStrEst} to hold for $\alpha > \frac{2}{q} + \frac{n}{p} - \frac{n}{2}$ and for $n>2$ we find the estimate \eqref{eq:homRegStrEst} to hold for $\alpha = \frac{2}{q} + \frac{n}{p}- \frac{n}{2}$, whenever $q \neq 2$.
\end{theorem}
 \section{Proof of Proposition \ref{prop:sharpRegularity} and Theorem \ref{thm:sharpRegularity}}
First, we show the necessary angular regularity:
\begin{proof}[{Proof of Proposition \ref{prop:sharpRegularity}}]
To simplify the concrete computation let us consider the special case of Schr\"odinger's equation, that is $\varphi(\rho) = \rho^2$ in \eqref{eq:homogeneousEquation}. Later on, we shall see that the example allows broader application.\\
As initial datum we consider a small rectangular block in frequency space
\begin{equation}
\label{eq:knappblock}
\hat{u}_0(\xi) = \chi_{(1-\varepsilon, 1+\varepsilon)}(\xi_1) \chi_{(-\varepsilon, \varepsilon)}(\xi_2) \ldots \chi_{(-\varepsilon, \varepsilon)}(\xi_n).
\end{equation}
For the solution to \eqref{eq:homogeneousEquation} we obtain
\begin{equation}
u(t,x) = C_n \int_{1-\varepsilon}^{1+\varepsilon} e^{i (x_1 \xi_1 - t \xi_1^2)} d \xi_1 \int_{-\varepsilon}^{\varepsilon} e^{i (x_2 \xi_2 - t \xi_2^2)} d \xi_2 \ldots \int_{-\varepsilon}^{\varepsilon} e^{i (x_n \xi_n - t \xi_n^2)} d \xi_n.
\end{equation}
We perform a change of variables $\tilde{\xi_1}  = \xi_1 -1$ and the first integral becomes
\begin{equation}
\int_{-\varepsilon}^{\varepsilon} d \tilde{\xi}_1 e^{i x_1 (1+ \tilde{\xi}_1) - i t ( \tilde{\xi}_1 +1)^2} = e^{i x_1} e^{- it} \int_{- \varepsilon}^{\varepsilon} d \tilde{\xi}_1 e^{i (x_1 -2t) \tilde{\xi}_1} e^{-it \tilde{\xi}_1^2 }.
\end{equation}
Observe that the phase nearly vanishes in the domain of integration, if $|x_1-2t| \ll 1/\varepsilon,\; |x_i| \ll 1/\varepsilon \,(i=2,\ldots,n), \; |t| \ll 1/\varepsilon^2$ and therefore we have got\footnote{This region is determined by the uncertainty principle $(\Delta x) (\Delta p) \gtrsim 1$ and by the group velocity $\varphi^\prime(1)$, where $\varphi$ denotes the dispersion relation.}
\begin{equation}
\label{eq:knappregion}
|u(t,x)| \gtrsim \varepsilon^{n} \mbox{ for} \; |x_1-2t| \ll 1/\varepsilon,\; |x_i| \ll 1/\varepsilon \,(i=2,\ldots,n), \; |t| \ll 1/\varepsilon^2
\end{equation}
and we obtain by the same means of the computation in the last section
\begin{equation}
\Vert u \Vert_{L_t^q L_x^p} \gtrsim \varepsilon^{n - n/p - 2/q}.
\end{equation}
Observe that smoothing out the sharp transitions does not change the behaviour found above. For the following part, we consider the sharp borders to be mollified on a scale of $\varepsilon/10$. The mollified variations will be denoted by $\tilde{\chi}$. We denote the non-negative derivative of the transition region centered at the origin by $\psi_\varepsilon(x)$, which is symmetric and satisfies $\psi_\varepsilon(0) \sim \varepsilon^{-1}, \; supp(\psi_\varepsilon) \subseteq (-\varepsilon/10, \varepsilon/10), \; \int \psi_\varepsilon = 1$; the family $(\psi_\varepsilon)_{\varepsilon} $ is in fact an approximate identity.\\
To estimate the Sobolev norm of the initial data with incorporated angular regularity, we note that we can equivalently consider the Fourier transform of the initial data (cf. ) and that the norm of the radial part of the initial datum is estimated by $\Vert \hat{u}_0^{rad} \Vert_{L_x^2} \sim \varepsilon^{n-1/2}$ because the radial part is roughly supported in $r = (1-\varepsilon, 1+\varepsilon)$ with a norm of order $\varepsilon^{n-1}$, which comes from intersecting the approximate block with spheres, when the intersection is approximately a square of sidelength $\varepsilon$ due to $\varepsilon$ very small. We still have to compute $\Vert (-\Delta_{\omega})^{1/2} \hat{u}_0 \Vert_{L^2}^2 = \sum_{1\leq i < j \leq n} \Vert \Omega_{ij} \hat{u}_0 \Vert_{L^2}^2 $.\\
Let us take a look at the concrete pair $i=1,\, j=2$, the other cases can be handled similarly. For this computation we will not keep track of the variables which remain untouched by the concrete vector field.
\begin{equation*}
\begin{split}
&\Omega_{12} \tilde{\chi}_{(1-\varepsilon,1+\varepsilon)}(x_1) \tilde{\chi}_{(-\varepsilon, \varepsilon)}(x_2) \\
&= x_1 \tilde{\chi}_{(1-\varepsilon,1+\varepsilon)}(x_1) \left( \psi_\varepsilon(x_2 + \varepsilon) - \psi_\varepsilon(x_2-\varepsilon) \right) \\
 &- \left( \psi_\varepsilon(x_1 - (1-\varepsilon)) - \psi_\varepsilon(x_1 - (1+\varepsilon)) \right) x_2 \tilde{\chi}_{(-\varepsilon, \varepsilon)}(x_2) 
\end{split}
\end{equation*}
Taking the square we find three terms: For the first one, we have
\begin{equation}
x_1^2 \tilde{\chi}^2_{(1-\varepsilon,1+\varepsilon)}(x_1) \left( \psi_\varepsilon^2(x_2+\varepsilon)+\psi_\varepsilon^2(x_2-\varepsilon) \right).
\end{equation}
Observe that this is the term of leading order after integration, which gives a quantity $\sim \varepsilon^{n-2}$. The other contributions are of order $\sim \varepsilon^{n-1}$.\\
Likewise we find that the contribution from the other vector fields of the kind $\Omega_{1j}$ is also of order $\sim \varepsilon^{n-1}$, contributions from vector fields of the kind $\Omega_{ij}, \, 1 \neq i < j$ are of higher order. We find the asymptotic
\begin{equation*}
\varepsilon^{n - n/p - 2/q} \lesssim \varepsilon^{n/2 - 1}
\end{equation*}
to hold for any $(q,p)$ radially admissible as $\varepsilon \rightarrow 0$. But, we find that after adding angular regularity of $\alpha<1$ we have got $\Vert u_0 \Vert_{\dot{H}_\omega^{0,\alpha}} \lesssim \varepsilon^{n/2 - \alpha}$ by means of interpolation and we find as a necessary condition
\begin{equation*}
\varepsilon^{n - n/p - 2/q} \lesssim \varepsilon^{n/2 - \alpha},
\end{equation*}
which gives as $\varepsilon \rightarrow 0$ that $\alpha \geq \frac{n}{p} + \frac{2}{q} - \frac{n}{2}$.
\end{proof}
Closely linked to estimates of the kind \eqref{eq:homRegStrEst} are estimates involving spherical averages, that is we consider the following norms:
 \begin{equation}
 \label{eq:sphericallyAveragedSolutions}
 \Vert u \Vert_{L_t^q \mathcal{L}_r^p L_\omega^2} = \left( \int_{\mathbb{R}} dt  \left( \int_0^\infty dr r^{n-1} \left( \int_{\mathbb{S}^{n-1}} d\omega |u(t,r \omega)|^2 \right)^{p/2} \right)^{q/p} \right)^{1/q}
 \end{equation}
The link is established through Sobolev embedding on the sphere 
 \begin{equation}
 \label{eq:SobolevEmbeddingSphere}
 \Vert f \Vert_{L_\omega^q} \lesssim_{n,q,\alpha} \Vert \Lambda_\omega^{\alpha} f \Vert_{L_\omega^2} \; \; ( q < \infty, \frac{1}{q} \geq \frac{1}{2} - \frac{\alpha}{n-1})
 \end{equation}
 and for Schr\"odinger-like equations Guo showed that those estimates hold in essentially the same range as the estimates found after requiring additional angular regularity:
  \begin{theorem}[{\cite[Theorem~1.1,~p.~3]{Guo2014}}]
 \label{thm:sphAvgStrichartzEstimates}
 Let $a > 1, \; q,p \geq 2$ and suppose that 
 \begin{equation*}
 \begin{split}
 &\frac{n}{2} \left( \frac{1}{2} - \frac{1}{p} \right) < \frac{1}{q} < \frac{2n-1}{2} \left( \frac{1}{2} - \frac{1}{p} \right) \mbox{ for } n =2, \\
 &\frac{n}{2} \left( \frac{1}{2} - \frac{1}{p} \right) < \frac{1}{q} \leq \frac{2n-1}{2} \left( \frac{1}{2} - \frac{1}{p} \right), \; (q,p) \neq \left( 2,\frac{4n-2}{2n-3} \right) \mbox{ for } n > 2.
 \end{split}
 \end{equation*}
Then we find the estimate 
\begin{equation}
\label{eq:sphAvgStrEst}
 \left\Vert P_N e^{it D^a} u_0 \right\Vert_{L_t^q \mathcal{L}_r^p L_\omega^2} \lesssim_{n,p,q} N^{s} \Vert u_0 \Vert_{L^2}
\end{equation} 
 to hold for any $N \in 2^{\mathbb{Z}}$ with $s =  \frac{n}{2} - \frac{n}{p} - \frac{a}{q} $.
 \end{theorem}
Theorem \ref{thm:sharpRegularity} is a consequence:
\begin{proof}[{Proof of Theorem \ref{thm:sharpRegularity}}]
For $(q,p)$ admissible for Theorem \ref{thm:sphAvgStrichartzEstimates} and $\alpha = (n-1) \left( \frac{1}{2} - \frac{1}{p} \right)$ we find
\begin{equation*}
\begin{split}
\Vert P_N u \Vert_{L_t^q L_x^p} &\lesssim_{n,p} \Vert \Lambda_\omega^\alpha P_N u \Vert_{L_t^q \mathcal{L}_r^p L_\omega^2} \\
&\lesssim_{n,p,q} N^{s} \Vert \Lambda_\omega^\alpha P_N u_0 \Vert_{L^2}
\end{split}
\end{equation*}
Taking $p$ to the sharp line, we find for $q>2$
\begin{equation*}
\alpha = (n-1) \left( \frac{1}{2} - \frac{1}{p} \right) = \frac{2}{q} + \frac{n}{p} - \frac{n}{2}.
\end{equation*}
For $q=2$ we only find
\begin{equation*}
\alpha > \frac{2}{q} + \frac{n}{p} - \frac{n}{2} \mbox{ as } p \rightarrow \frac{4n-2}{2n-3}.
\end{equation*}
The frequency localization is removed employing Littlewood-Paley theory and the proof is concluded interpolating with the estimates on the classical sharp line.
\end{proof}
\section{Remarks}
 The special case of the wave equation, which corresponds to $a=1$ in \eqref{eq:homogeneousEquation}, had been treated in \cite{Cho2013,Sterbenz2005} and the results holding the analogues of the estimates we discussed above are sharp up to endpoints. In \cite[Corollary~1.5.,~p.~4]{ChoGuoLee2015} the endpoint angular regularity was obtained for $q \geq p$.\\
From the weakened dispersion the classical estimates, which were already found in \cite{TaoKeel1998}, are only valid in the range
 \begin{equation*}
 \frac{1}{q} \leq \frac{n-1}{2}  \left( \frac{1}{2} - \frac{1}{p} \right), \; q,p \geq 2, \; p \neq \infty,
 \end{equation*}
 whereat the estimates found after taking spherical averages or adding angular regularity are valid in the range
 \begin{equation*}
 \frac{n-1}{2}  \left( \frac{1}{2} - \frac{1}{p} \right) < \frac{1}{q} < (n-1) \left( \frac{1}{2} - \frac{1}{p} \right), \; q,p \geq 2, \; p \neq \infty.
 \end{equation*}
Sharpness is again found employing Knapp-type examples (cf. \cite[pp.~190,~195]{Sterbenz2005}). In \cite{Sterbenz2005} was also considered additional angular regularity.\\
Furthermore, we note that Proposition \ref{prop:sharpRegularity} also provides a Knapp-type example for far more general dispersion relations than those associated to Schr\"odinger-like equations, for which Cho and Lee proved estimates of the kind \eqref{eq:homRegStrEst} in the same range like in Theorem \ref{thm:radStrEstA} provided that $\alpha > \frac{5n-1}{5n-5} \left( \frac{2}{q} + \frac{n}{p} - \frac{n}{2} \right)$ (cf. \cite[Theorem~1.1.,~p.~993]{Cho2013}).\\
The example given in Proposition \ref{prop:sharpRegularity} leads us to conjecture that even for more general dispersion relations the sharp angular regularity is $\alpha = \frac{2}{q} + \frac{n}{p} - \frac{n}{2}$.

\bibliographystyle{amsxport}
\end{document}